\documentclass[12pt]{amsart}
\usepackage{amscd}
\usepackage{amsmath}
\usepackage{amssymb}
\usepackage{amsthm}
\newtheorem{thm}{Theorem}[section]

\newtheorem{lemma}[thm]{Lemma}

\newtheorem{cor}[thm]{Corollary}

\theoremstyle{definition}
\newtheorem{defn}[thm]{Definition}
\newtheorem{exam}[thm]{Example}

\newtheorem*{acknowledgement}{Acknowledgments}
\theoremstyle{remark}
\newtheorem{remark}[thm]{Remark}
\numberwithin{equation}{section}

\DeclareMathOperator{\pd}{pd}

\DeclareMathOperator{\bight}{bight}

\DeclareMathOperator{\reg}{reg}

\DeclareMathOperator{\lcm}{lcm}

\DeclareMathOperator{\dist}{dist}

\begin{document}
\title[Non-vanishingness of Betti numbers]
  {Non-vanishingness of Betti numbers of edge ideals}
\author[Kyouko Kimura]{Kyouko Kimura}
\address
         {Department of Mathematics, Faculty of Science, 
          Shizuoka University, 836 Ohya, Suruga-ku, Shizuoka 422-8529, Japan}
\email{skkimur@ipc.shizuoka.ac.jp}
\date{}
\begin{abstract}
  Given finite simple graph one can associate the edge ideal. 
  In this paper we discuss the non-vanishingness of the graded Betti numbers 
  of edge ideals in terms of the original graph. 
  In particular, we give a necessary and sufficient condition 
  for a chordal graph on which the graded Betti number does not vanish 
  and characterize the graded Betti number for a forest. 
  Moreover we characterize the projective dimension for a chordal graph. 
\end{abstract}
\maketitle

\section{Introduction}
Let $G$ be a finite simple graph, 
that is, finite graph with no loop and no multiple edge. 
We denote its vertex set by $V=V(G)$ and its edge set by $E(G)$. 
Let $K$ be a field and $S=K[V]:=K[x : x \in V]$ a polynomial ring 
with $\deg x = 1$. 
The \textit{edge ideal} of $G$ is a squarefree monomial ideal $I(G) \subset S$ 
generated by all products $x_i x_j$ with $\{ x_i, x_j \} \in E(G)$. 
We are interested in describing invariants of $I(G)$ in terms of $G$. 

\par
Let us consider a minimal graded free resolution of $S/I(G)$ over $S$: 
\begin{displaymath}
  0 \longrightarrow \bigoplus_{j} S(-j)^{\beta_{p,j}} 
    \longrightarrow \cdots
    \longrightarrow \bigoplus_{j} S(-j)^{\beta_{1,j}}
    \longrightarrow S \longrightarrow S/I(G) \longrightarrow 0. 
\end{displaymath}
The integers $\beta_{i,j} (S/I(G)) := \beta_{i,j}$ is called
the $i$th \textit{graded Betti number} of $S/I(G)$ in degree $j$. 
The length $p$ of the resolution is called the projective dimension of 
$S/I(G)$ over $S$, denoted by $\pd S/I(G)$, that is, 
\begin{displaymath}
  \pd S/I(G) 
   = \max \{ i : \beta_{i,j} (S/I(G)) \neq 0 \  \text{for some $j$} \}. 
\end{displaymath}
Also, the \textit{(Castelnuovo--Mumford) regularity} of $S/I(G)$ is defined by 
\begin{displaymath}
  \reg S/I(G) := \max \{ j-i : {\beta}_{i,j} (S/I(G)) \neq 0 \}. 
\end{displaymath}
In this paper, we focus on these invariants, 
which have studied by many authors, 
e.g., \cite{Barile08, HaTuyl08, Jacques, JacqKatz, 
Katzman, Kummini, MMCRTY, Zheng}. 
In particular, 
Zheng \cite{Zheng} characterized the projective dimension 
and the regularity for a \textit{forest}, which is a graph with no cycle. 
Later H\`{a} and Van Tuyl \cite{HaTuyl08} extended 
this characterization of the regularity to that for a chordal graph. 
Here a finite graph is called \textit{chordal} 
if each cycle of $G$ whose length is more than $3$ has a chord. 
On the other hand, Katzman \cite{Katzman} proved some results on 
the non-vanishingness of the graded Betti numbers. 
For other problems and results on this area we refer to \cite{HaTuyl07}. 

\par
In this paper, we give a sufficient condition on which the graded Betti number 
does not vanish (Theorem \ref{claim:NonVanishBetti}). 
This is a generalization of Katzman's result. 
Moreover we proved that the condition is also a necessary one 
for chordal graphs (Theorem \ref{claim:Betti_chordal}). 
As a result, we characterize the non-vanishingness of the graded Betti numbers 
and thus, the projective dimension for chordal graphs. 
Furthermore we give a characterization of the graded Betti numbers for forests. 

\par
The organization of the paper is as follows. 
In Section 2, we recall some definitions on graphs and results 
by Zheng \cite{Zheng}, H\`{a} and Van Tuyl \cite{HaTuyl08}, 
and Katzman \cite{Katzman} mentioned above. 
Also in this section we introduce the notion of the 
\textit{strongly disjoint set of bouquets} on a graph 
which plays an important role in our characterization. 
Next in Section 3, we discuss the non-vanishingness 
of the graded Betti numbers. 
In particular, we prove Theorem \ref{claim:NonVanishBetti}. 
In Section 4, we provide refined results 
of Theorem \ref{claim:NonVanishBetti} for chordal graphs or forests: 
Theorem \ref{claim:Betti_chordal}. 
Then we also have a characterization of the projective dimension 
for a chordal graph. 
However, in Section 5, we give another characterization of it 
which has a rather weaker condition. 
As an application of this, we recover the result that 
the projective dimension of the edge ideal of a chordal graph coincides 
with its big height (Corollary \ref{claim:bight}), 
which was proved by Morey and Villarreal \cite[Corollary 3.33]{MV} 
for more general graphs. 

\section{Strongly disjoint set of bouquets}
In this section, we prepare some definitions on graphs 
and recall some known results on our problem. 

\par
Let $G$ be a finite simple graph with the vertex set $V$ 
and the edge set $E(G)$. 
Let $e, e'$ be two distinct edges of $G$. 
The \textit{distance} between $e$ and $e'$ in $G$, 
denoted by $\dist_G (e,e')$, 
is defined by the minimum length $\ell$ among sequences 
$e_0 = e, e_1, \ldots, e_{\ell} = e'$ 
with $e_{i-1} \cap e_i \neq \emptyset$, 
where $e_i \in E(G)$. If there is no such a sequence, we define 
$\dist_G (e, e') = \infty$. 
We say that $e$ and $e'$ are \textit{$3$-disjoint} in $G$ 
if $\dist_G (e, e') \geq 3$. 
A subset $\mathcal{E} \subset E(G)$ is said to be 
\textit{pairwise $3$-disjoint} if every pairs of distinct edges 
$e, e' \in \mathcal{E}$ are $3$-disjoint in $G$; 
see \cite[Definitions 2.2 and 6.3]{HaTuyl08}. 

\par
The graph $B$ with $V(B) = \{w, z_1, \ldots, z_d \}$ 
and $E(B) = \{ \{ w, z_i \} : i = 1, \ldots, d \}$ ($d \geq 1$)
is called a \textit{bouquet}. 
  
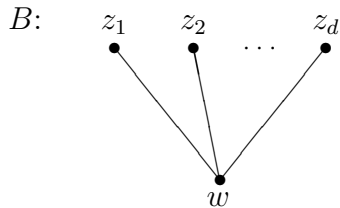
\begin{center}
    \begin{picture}(130,80)(-30,10)
      \put(-30,77){$B$:}
      \put(50,20){\circle*{4}}
      \put(45,10){$w$}
      \put(10,70){\circle*{4}}
      \put(5,77){$z_1$}
      \put(40,70){\circle*{4}}
      \put(35,77){$z_2$}
      \put(60,70){\circle*{1}}
      \put(65,70){\circle*{1}}
      \put(70,70){\circle*{1}}
      \put(90,70){\circle*{4}}
      \put(85,77){$z_d$}
      \put(50,20){\line(-4,5){40}}
      \put(50,20){\line(-1,5){10}}
      \put(50,20){\line(4,5){40}}
    \end{picture}
  \end{center}
Then the vertex $w$ is called the \textit{root} of $B$, the vertices $z_i$ 
\textit{flowers} of $B$, and the edges $\{ w, z_i \}$ \textit{stems} of $B$; 
see \cite[Definition 1.7]{Zheng}. 
We call a bouquet which is a subgraph of 
a graph $G$ a bouquet of $G$. 
Let $\mathcal{B} = \{ B_1, B_2, \ldots, B_j \}$ be a set of bouquets of $G$. 
We set 
\begin{displaymath}
  \begin{aligned}
    F({\mathcal{B}}) &:= \{ z \in V \; : \; 
      \text{$z$ is a flower of some bouquet in $\mathcal{B}$} \}, \\
    R({\mathcal{B}}) &:= \{ w \in V \; : \; 
      \text{$w$ is a root of some bouquet in $\mathcal{B}$} \}, \\
    S({\mathcal{B}}) &:= \{ s \in E(G) \; : \; 
      \text{$s$ is a stem of some bouquet in $\mathcal{B}$} \}. 
  \end{aligned}
\end{displaymath}
A \textit{type} of $\mathcal{B}$ is defined by 
$(\# F({\mathcal{B}}), \# R({\mathcal{B}}))$. 
We define a disjointness on the set of bouquets of $G$. 
\begin{defn}
  \label{defn:StrongDisjoint}
  A set $\mathcal{B} = \{ B_1, B_2, \ldots, B_j \}$ of bouquets of $G$ 
  is said to be \textit{strongly disjoint} in $G$ if the following 
  $2$ conditions are satisfied: 
  \begin{enumerate}
  \item $V(B_k) \cap V(B_{\ell}) = \emptyset$ for all $k \neq \ell$. 
  \item We can choose a stem $s_k$ from each bouquets $B_k \in \mathcal{B}$ 
    so that $\{ s_1, s_2, \ldots, s_j \}$ is pairwise $3$-disjoint in $G$. 
  \end{enumerate}
\end{defn}
%
%
\begin{remark}
  \label{rmk:root}
  If $\mathcal{B} = \{ B_1, B_2, \ldots, B_j \}$ is a strongly disjoint set 
  of bouquets in $G$, then any two vertices those belong to 
  $R({\mathcal{B}})$ is not adjacent in $G$. 
  Indeed if not, say the roots of $B_1$ and $B_2$ are adjacent in $G$, 
  then the distance between any stem of $B_1$ and any stem of $B_2$ is $2$. 
\end{remark}

\par
Moreover we give the following definition. 
\begin{defn}
  \label{defn:ContainStrongDisjoint}
  Let $G$ be a finite simple graph. 
  \begin{enumerate}
  \item We say that \textit{$G$ coincides with a strongly disjoint set 
    of bouquets of type $(i,j)$} if there exists a strongly disjoint set 
    $\mathcal{B}$ of bouquets of type $(i,j)$ with 
    \begin{displaymath}
      V(G) = F({\mathcal{B}}) \cup R ({\mathcal{B}}), \qquad 
      E(G) = S({\mathcal{B}}). 
    \end{displaymath}
  \item We say that \textit{$G$ contains a strongly disjoint set 
    of bouquets of type $(i,j)$} if there exists a strongly disjoint set 
    $\mathcal{B}$ of bouquets of type $(i,j)$ with 
    \begin{displaymath}
      V(G) = F({\mathcal{B}}) \cup R ({\mathcal{B}}). 
    \end{displaymath}
  \end{enumerate}
\end{defn}

\par
Let $G$ be a finite simple graph on $V$ and $W$ a subset of $V$. 
The subgraph of $G$ whose vertex set is $W$ and whose edge set is 
$\{ e \in E(G) : e \subset W \}$ is called the induced subgraph of $G$ on $W$ 
and denoted by $G_W$. In the forward sections we consider an induced subgraph 
of $G$ which contains a strongly disjoint set of bouquets 
of type $(i,j)$. 

\par
Now we recall the results due to Zheng \cite{Zheng} and 
H\`{a} and Van Tuyl \cite{HaTuyl08}, 
mentioned in Introduction. 
For a graph $G$, we set 
\begin{displaymath}
  \begin{aligned}
    d_G &:= \max \{ \# F({\mathcal{B}}) \; : \; 
      \text{$\mathcal{B}$ is a strongly disjoint set of bouquets of $G$} \}, \\
    c_G &:= \max \{ \# \mathcal{E} \; : \; 
      \text{$\mathcal{E} \subset E(G)$ is a pairwise $3$-disjoint in $G$} \}. 
  \end{aligned}
\end{displaymath}
\begin{thm}[{H\`{a} and Van Tuyl \cite{HaTuyl08}, Zheng \cite{Zheng}}]
  \label{claim:regularity}
  Let $G$ be a finite simple graph on $V$ and $S=K[V]$. 
  \begin{enumerate}
  \item $($\cite{Zheng}$)$. When $G$ is a forest, we have $\pd S/I(G) = d_G$. 
  \item $($\cite[Theorem 6.5]{HaTuyl08}$)$. $\reg S/I(G) \geq c_G$. 
  \item $($\cite[Theorem 6.8]{HaTuyl08}$,$ \cite[Theorem 2.18]{Zheng}$)$. 
    When $G$ is a chordal graph, 
    we have $\reg S/I(G) = c_G$. 
  \end{enumerate}
\end{thm}

\par
Also Katzman's results on graded Betti numbers are as follows: 
\begin{thm}[{Katzman \cite[Lemma 2.2, Proposition 2.5]{Katzman}}]
  \label{claim:Katzman}
  Let $G$ be a finite simple graph on $V$ and $S=K[V]$. 
  \begin{enumerate}
  \item If there exists a subset $W \subset V$ such that the induced subgraph 
    $G_W$ coincides with a strongly disjoint set of bouquets of type 
   $(i,j)$, then $\beta_{i,i+j} (S/I(G)) \neq 0$. 
  \item $\beta_{i,k} (S/I(G)) = 0$ when $k > 2i$. 
  \item The graded Betti number $\beta_{i,2i} (S/I(G))$ 
    coincides with the number of subsets $W$ of $V$ for which 
    the induced subgraph $G_W$ coincides with a strongly disjoint 
    set of bouquets of type $(i,2i)$. 
  \end{enumerate}
\end{thm}

\begin{remark}
  \label{rmk:hypergraph}
  Theorem \ref{claim:regularity} (2) follows from 
  Theorem \ref{claim:Katzman} (1) or (3). 
  Actually the result of H\`{a} and Van Tuyl \cite{HaTuyl08} is 
  in more general situation of hypergraphs. 
\end{remark}

\section{Non-vanishingness of the graded Betti numbers}
In this section we provide a sufficient condition 
for which a graded Betti number of an edge ideal does not vanish. 
The main result in this section is the following theorem, 
which is a generalization 
of Katzman's result (Theorem \ref{claim:Katzman} (1)). 
\begin{thm}
  \label{claim:NonVanishBetti}
  Let $G$ be a finite simple graph on $V$ and $S = K[V]$. 
  Assume that there exists a subset $W$ of $V$ such that 
  the induced subgraph $G_W$ contains a strongly disjoint set of 
  bouquets of type $(i,j)$. Then $\beta_{i, i+j} (S/I(G)) \neq 0$. 

  \par
  In particular, we have $\pd S/I(G) \geq d_G$. 
\end{thm}
\begin{remark}
  Precisely we have that $\beta_{i, i+j} (S/I(G))$ is greater than or equal to 
  the number of subsets $W \subset V$ those satisfy the condition in 
  Theorem \ref{claim:NonVanishBetti}. 
\end{remark}

\par
Before proving Theorem \ref{claim:NonVanishBetti}, 
we investigate a relation between the graded Betti numbers of the edge ideal of 
a graph and that of the induced subgraphs. 
Note that $\beta_{i,j} (K[V(G)]/I(G)) = \beta_{i-1,j} (I(G))$ for $i \geq 1$. 
\begin{lemma}
  \label{claim:restriction}
  Let $G$ be a finite simple graph on $V$. Then for all $i \geq 0$, we have 
  \begin{displaymath}
    \beta_{i,j} (I(G)) = \sum_{\genfrac{}{}{0pt}{}{W \subset V}{\# W = j}} 
      \beta_{i,j} (I(G_W)). 
  \end{displaymath}
\end{lemma}
\begin{proof}
  Since $I (G)$ is a squarefree monomial ideal, 
  it is the Stanley--Reisner ideal $I_{\Delta}$ for some simplicial complex 
  $\Delta$. 
  By Hochster's formula for Betti numbers 
  (see e.g., \cite[Theorem 5.5.1]{BrunsHerzog}), we have 
  \begin{displaymath}
    \beta_{i,j}(I(G)) = \beta_{i+1,j} (K[{\Delta}]) 
    = \sum_{\genfrac{}{}{0pt}{}{W \subset V}{\# W = j}}
      \dim_K \tilde{H}_{j-(i+1)-1}({\Delta}_W; K), 
  \end{displaymath}
  where $\tilde{H}_i({\Delta}; K)$ stands for the $i$th reduced homology 
  group of $\Delta$ and where $\Delta_W$ 
  denotes the restriction of $\Delta$ on $W$: 
  $\Delta_W = \{ F \in \Delta : F \subset W \}$. 
  Note that $I_{{\Delta}_W} = I(G_W)$ for a subset $W \subset V$. 
  Hence again by Hochster's formula, we have 
  \begin{displaymath}
    \beta_{i,j}(I(G_W)) = \beta_{i+1,j} (K[{\Delta}_W]) 
      = \dim_K \tilde{H}_{j-(i+1)-1}({\Delta}_{W}; K) 
  \end{displaymath}
  for a subset $W \subset V$ with $\# W = j$. 
   This completes the proof. 
\end{proof}

\par
In the proof of Theorem \ref{claim:Katzman} (1), 
Katzman used the Taylor resolution. 
In the proof of Theorem \ref{claim:NonVanishBetti}, 
we use a Lyubeznik resolution in stead of the Taylor resolution. 
A Lyubeznik resolution (\cite{Ly}) is a subcomplex of the Taylor resolution 
generated by $L$-admissible symbols. It gives a 
(not necessarily minimal) free resolution for a monomial ideal. 

\par
Let $I$ be a monomial ideal of a polynomial ring over $K$ and 
$m_1, m_2, \ldots, m_{\mu}$ the minimal system of monomial generators of $I$. 
Let $e_{{\ell}_1 {\ell}_2 \cdots {\ell}_i}$ 
(${\ell}_1 < {\ell}_2 < \cdots < {\ell}_i$) denotes 
the free basis of the Taylor resolution of $I$. 
Recall that the degree of $e_{{\ell}_1 {\ell}_2 \cdots {\ell}_i}$ is given by 
the degree of $\lcm (m_{{\ell}_1}, m_{{\ell}_2}, \ldots, m_{{\ell}_i})$. 
We say that a symbol $[{\ell}_1, {\ell}_2, \ldots, {\ell}_i] 
  := e_{{\ell}_1 {\ell}_2 \cdots {\ell}_i}$ is 
\textit{$L$-admissible} if for all $t<i$ and for all $q < {\ell}_t$, 
the monomial generator $m_q$ does not divide 
$\lcm (m_{{\ell}_t}, m_{{\ell}_{t+1}}, \ldots, m_{{\ell}_i})$. 
Note that the $L$-admissibleness as well as a Lyubeznik resolution of $I$ 
depends on an order of monomial generators of $I$. 
An $L$-admissible symbol $[{\ell}_1, {\ell}_2, \ldots, {\ell}_i]$ 
is said to be maximal if there is no $L$-admissible symbol 
$[k_1, k_2, \ldots, k_{i'}]$ with 
$\{ {\ell}_1, {\ell}_2, \ldots, {\ell}_i \} \subsetneq \{ k_1, k_2, \ldots, k_{i'} \}$. 
If there exists a maximal $L$-admissible symbol 
$[{\ell}_1, {\ell}_2, \ldots, {\ell}_i]$ which satisfies the condition
\begin{equation}
  \label{eq:maxLadmissible}
  \lcm (m_{{\ell}_1}, \ldots, \widehat{m_{{\ell}_k}}, \ldots, m_{{\ell}_i}) 
  \neq \lcm (m_{{\ell}_1}, \ldots, m_{{\ell}_i}), 
  \quad \text{for all $k=1, 2, \ldots, i$}, 
\end{equation}
then $\beta_{i,j} (S/I(G)) \neq 0$, where
$j=\deg [{\ell}_1, {\ell}_2, \ldots, {\ell}_i]$ 
(see also Barile \cite[Remark 1]{Barile05}). 

\par
\bigskip

\par
Now we prove Theorem \ref{claim:NonVanishBetti}. 
\begin{proof}[Proof of Theorem \ref{claim:NonVanishBetti}]
  By virtue of Lemma \ref{claim:restriction}, 
  it is sufficient to prove the theorem 
  when $\# V = i+j$ and $G$ contains a strongly disjoint set of 
  bouquets of type $(i,j)$. Let $\mathcal{B} = \{ B_1, B_2, \ldots, B_j \}$ 
  be such a set of bouquets. Since $\mathcal{B}$ is strongly disjoint, 
  we can choose a stem $s_k$ from each bouquet $B_k \in \mathcal{B}$ 
  so that ${\mathcal{S}}_0 := \{ s_1, s_2, \ldots, s_j \}$ is 
  pairwise $3$-disjoint in $G$. 
  We set $\mathcal{S}' = S({\mathcal{B}}) \setminus {\mathcal{S}}_0$ 
  and $\mathcal{E} = E(G) \setminus S({\mathcal{B}})$. 
  We define an ordering of the edges of $G$ as 
  $\mathcal{S}', \mathcal{E}, {\mathcal{S}}_0$ and consider the associated 
  ordering on the minimal system of monomial generators of $I(G)$. 
  We consider the symbol $\sigma$ corresponding to 
  $\mathcal{S}', \mathcal{S}_0$. We claim that $\sigma$ 
  is a maximal $L$-admissible symbol. 

  \par
  The $L$-admissibleness of $\sigma$ follows from the assumption 
  that $\mathcal{S}_0$ is pairwise $3$-disjoint. 
  To prove that $\sigma$ is maximal, we consider the symbol $\tau$ which 
  corresponds to $\mathcal{S}', e, \mathcal{S}_0$, where 
  $e = \{ u, v \} \in \mathcal{E}$. 
  Since $V = F({\mathcal{B}}) \cup R({\mathcal{B}})$ and 
  $\mathcal{B}$ is strongly disjoint, it follows that 
  $\{ u, v \} \cap F({\mathcal{B}}) \neq \emptyset$. 
  Moreover at least one of $u, v$ belongs to $F({\mathcal{B}})$ which is not 
  a vertex of the stems belonging to ${\mathcal{S}}_0$ 
  because ${\mathcal{S}}_0$ is pairwise $3$-disjoint. 
  Let $u$ be such a vertex and assume that $u \in V(B_k)$. 
  Then the product of the monomial $uv$ and the monomial which corresponds to 
  $s_k$ is divisible by the monomial corresponding to the stem of $B_k$ 
  whose flower is $u$. Hence $\tau$ is not $L$-admissible. 
  Therefore $\sigma$ is a maximal $L$-admissible symbol. 

  \par
  It is clear that $\sigma$ satisfies the condition (\ref{eq:maxLadmissible}). 
  Since $\deg \sigma = i+j$, we conclude that $\beta_{i, i+j} (S/I(G)) \neq 0$. 
\end{proof}

\section{The case of chordal graphs}
In this section, we prove that the converse 
of Theorem \ref{claim:NonVanishBetti} is true when $G$ is a chordal graph. 
Moreover we give a characterization of the graded Betti numbers 
of edge ideals of forests. 
Precisely, the main result of this section is the following theorem. 
\begin{thm}
  \label{claim:Betti_chordal}
  Let $G$ be a finite simple graph on $V$ and $S=K[V]$. 
  \begin{enumerate}
  \item Suppose that $G$ is chordal. Then $\beta_{i, i+j} (S/I(G)) \neq 0$ 
    if and only if there exists a subset $W$ of $V$ such that 
    the induced subgraph $G_W$ contains a strongly disjoint set of 
    bouquets of type $(i,j)$. 

    \par
    In particular, we have $\pd S/I(G) = d_G$. 
  \item When $G$ is a forest, the graded Betti number $\beta_{i, i+j} (S/I(G))$ 
    coincides with the number of subsets $W$ of $V$ with the same condition 
    as in (1). 
  \end{enumerate}
\end{thm}
\begin{remark}
  We obtain a characterization of the projective dimension 
  for chordal graphs by Theorem \ref{claim:Betti_chordal} 
  though we give another characterization in the next section. 
\end{remark}

\par
Let $G$ be a finite simple graph on $V$. Take a vertex $v \in V$. 
We say that $u \in V$ is a neighbour of $v$ in $G$ if 
$\{ u,v \}$ is an edge of $G$. 
We denote by $N(v)$, the neighbourhood of $v$. 
For an edge  $e \in E(G)$, we denote by $G \setminus e$ the subgraph of $G$ 
obtained from $G$ by removing the edge $e$. 

\par
A graph $G$ is said to be \textit{chordal} if each cycle in $G$ 
whose length is more than $3$ has a chord. 
Dirac \cite{Dirac} proved that when $G$ is chordal, 
there exists a perfect elimination ordering 
on $E(G)$. This means that any induced subgraph $G'$ of $G$ has a vertex $v$ 
such that the induced subgraph of $G'$ 
on the neighbourhood of $v$ in $G'$ 
is a complete graph. 

\par
The next result about the graded Betti numbers of edge ideals of chordal graphs 
due to H\`{a} and Van Tuyl 
is a key in our proof of Theorem \ref{claim:Betti_chordal}. 
For simplicity, we set $\beta_{-1,0} (I(G)) = 1$ 
and $\beta_{-1, j} (I(G)) = 0$ if $j \neq 0$. 
When $G$ is a graph with no edge, we set $\beta_{i,j} (I(G))$ as $1$ when 
$(i,j) = (-1,0)$ and as $0$ otherwise. 
\begin{lemma}[{H\`{a} and Van Tuyl \cite[Theorem 5.8]{HaTuyl08}}]
  \label{claim:HaTuyl}
  Let $G$ be a chordal graph on the vertex set $V$. 
  Suppose that $e = \{ u,v \}$ is an edge 
  of $G$ such that $G_{N(v)}$ is a complete graph. 
  Set $N(u) = \{ v, x_1, \ldots, x_t \}$ and 
  $G' = G_{V \setminus \{ u, v, x_1, \ldots, x_t \}}$. 
  Then both $G \setminus e$ and $G'$ are chordal, and 
  \begin{equation}
    \label{eq:HaTuyl}
    \beta_{i,j} (I(G)) = \beta_{i,j}(I(G\setminus e)) 
      + \sum_{\ell = 0}^{i} \binom{t}{\ell} 
        \beta_{i - 1 - \ell, j - 2 - \ell} (I(G')). 
  \end{equation}
\end{lemma}

\begin{remark}
  \label{rmk:G'}
  The edge set of $G'$ in the above theorem is 
  \begin{displaymath}
    E(G') = \{ e' \in E(G) \; : \; \dist_G (e, e') \geq 3 \}. 
  \end{displaymath}
\end{remark}

\par
In the proof of Theorem \ref{claim:Betti_chordal}, we use (\ref{eq:HaTuyl}) 
as the following form: 
  \begin{equation}
    \label{eq:HaTuylproof}
    \beta_{i-1,i+j} (I(G)) = \beta_{i-1,i+j}(I(G\setminus e)) 
      + \sum_{\ell = 0}^{i-1} \binom{t}{\ell} 
        \beta_{i - 2 - \ell, (i - 1 - \ell) + (j-1)} (I(G')). 
  \end{equation}

\par
First we investigate the relation between a strongly disjoint set of bouquets
of $G \setminus e, G'$ and that of $G$.  
\begin{lemma}
  \label{claim:d}
  Let $G$ be a chordal graph. 
  We use the same notations as in Lemma \ref{claim:HaTuyl}. 
  \begin{enumerate}
  \item Let $\mathcal{B}$ be a strongly disjoint set of bouquets of 
    $G \setminus e$. Then $\mathcal{B}$ is also strongly disjoint in $G$. 
    In particular, $d_G \geq d_{G \setminus e}$. 
  \item Let $\mathcal{B}'$ be a strongly disjoint set of bouquets of $G'$ 
    and $B$ the bouquet of $G$ whose root is $u$ 
    and whose flowers are $v, x_1, \ldots, x_t$. 
    Then $\mathcal{B}' \cup \{ B \}$ is a strongly disjoint set of bouquets 
    of $G$. 
    In particular, $d_G \geq d_{G'} + (t+1)$. 
  \end{enumerate}
\end{lemma}
\begin{proof}
  (1) Let $\mathcal{B} = \{ B_1, B_2, \ldots, B_j \}$ be 
  a strongly disjoint set of bouquets in $G \setminus e$, 
  where $e = \{ u, v \}$. 
  If one of $u,v$ is not in $R({\mathcal{B}}) \cup F({\mathcal{B}})$, 
  then it is clear that $\mathcal{B}$ is also strongly disjoint in $G$. 
  Also we easily see that $\mathcal{B}$ is strongly disjoint in $G$ when 
  $u, v$ are the vertices of the same bouquet $B_k$. (In this case, 
  both of $u, v$ are flowers of $B_k$.) 
  Thus we may assume that $u \in V(B_1)$ and $v \in V(B_2)$. 

  \par
  We consider $4$ cases. 

  \par 
  \textbf{Case 1:} $u, v \in R(\mathcal{B})$. 
  In this case, all flowers of $B_2$ are neighbours of $u$. 
  Then there are no stems $s_1, s_2$ ($s_1 \in E(B_1)$, $s_2 \in E(B_2)$) 
  with $\dist_{G \setminus e} (s_1, s_2) \geq 3$. 
  This contradicts the assumption that $\mathcal{B}$ is strongly disjoint in 
  $G \setminus e$. 

  \par 
  \textbf{Case 2:} $u \in R(\mathcal{B})$ and $v \in F(\mathcal{B})$. 
  Since $G_{N(v)}$ is a complete graph, the root of $B_2$ is a neighbour 
  of $u$. This leads a contradiction as in Case 1. 

  \par 
  \textbf{Case 3:} $u \in F(\mathcal{B})$ and $v \in R(\mathcal{B})$. 
  Let $w_1$ be the root of $B_1$. The completeness of 
  $G_{N(v)}$ implies that 
  the distance between $\{ u, w_1 \}$ and any stem $s_2$ of $B_2$ in 
  $G \setminus e$ is equal to $2$. 
  It follows that $\mathcal{B}$ is also strongly disjoint in $G$. 

  \par 
  \textbf{Case 4:} $u,v \in F(\mathcal{B})$. 
  Let $w_1, w_2$ be roots of $B_1, B_2$ respectively. 
  Then $\{ w_2, u \} \in E(G \setminus e)$ and 
  $\dist_{G \setminus e} (\{ u, w_1 \}, \{ v, w_2 \}) = 2$. 
  Therefore $\mathcal{B}$ is also strongly disjoint in $G$. 

  \par
  (2) We observed that $e$ is $3$-disjoint with any stem of $\mathcal{B}'$ 
  in Remark \ref{rmk:G'}. Let $\mathcal{S}$ be a set of stems of $\mathcal{B}'$ 
  which guarantees the strongly disjointness of $\mathcal{B}'$ in $G'$. 
  Then $\mathcal{S} \cup \{ e \}$ is pairwise $3$-disjoint. 
  This guarantees the strongly disjointness of $\mathcal{B}' \cup \{ B \}$. 
 \end{proof}

\par
Now we prove Theorem \ref{claim:Betti_chordal}. 
\begin{proof}[{Proof of Theorem \ref{claim:Betti_chordal}}]
  (1) By Theorem \ref{claim:NonVanishBetti}, it is sufficient to prove
  that when $\beta_{i, i+j} (S/I(G)) \neq 0$, 
  there exists a subset $W$ of $V$ such that $G_W$ contains 
  a strongly disjoint set of bouquets of type $(i,j)$. 
  We use induction on $\# E(G)$. 
  If $\# E(G) = 1$, then $\beta_{i, i+j} (S/I(G)) = 0$ 
  except for $(i,j) = (1,1)$ 
  and $\beta_{1, 1+1} (S/I(G)) = 1$. In this case, $G$ coincides with 
  a strongly disjoint set of bouquets of type $(1,1)$. Hence the claim is true. 

  \par
  Assume that $\# E(G) \geq 2$. 
  By Lemma \ref{claim:restriction}, 
  we may assume that $\# V = i+j$ and may prove that 
  if $\beta_{i,i+j} (S/I(G)) \neq 0$, i.e., $\beta_{i-1,i+j} (I(G)) \neq 0$, 
  then $G$ contains 
  a strongly disjoint set of bouquets of type $(i,j)$. 
  We use the same notation as in Lemma \ref{claim:HaTuyl}. 
  Note that $G' = G_{W_0}$ 
  where $W_0 = V \setminus \{ u, v, x_1, \ldots, x_t \}$. 
  Since $\# W_0 = i+j - (t+2)$, the summand of righthand side 
  of the second term of (\ref{eq:HaTuylproof}) is $0$ except for $\ell = t$. 
  Hence we can rewrite (\ref{eq:HaTuylproof}) as 
  \begin{equation}
    \label{eq:HaTuylproof2}
    \beta_{i-1,i+j} (I(G)) = \beta_{i-1,i+j}(I(G\setminus e))  
    + \beta_{i-2-t, (i-1-t)+(j-1)} (I(G')). 
  \end{equation}
  By the assumption that $\beta_{i-1, i+j} (I(G)) \neq 0$, 
  one of the summands of righthand side of (\ref{eq:HaTuylproof2}) 
  does not vanish. 

  \par
  If $\beta_{i-1, i+j} (I(G \setminus e)) \neq 0$, 
  then by inductive hypothesis, $G \setminus e$ contains 
  a strongly disjoint set of bouquets of type $(i,j)$. 
  By Lemma \ref{claim:d} (1), this set of bouquets 
  is also strongly disjoint in $G$. 

  \par
  If $\beta_{i-2-t, (i-1-t)+(j-1)} (I(G')) \neq 0$, 
  then by inductive hypothesis, $G'$ contains 
  a strongly disjoint set $\mathcal{B}'$ of bouquets of type $(i-1-t,j-1)$. 
  Let $B$ be the bouquet of $G$ whose root is $u$ and whose flowers are 
  $v, x_1, \ldots, x_t$. 
  By Lemma \ref{claim:d} (2), the set of bouquets 
  $\mathcal{B} := \mathcal{B}' \cup \{ B \}$
  is strongly disjoint in $G$ and the type is $(i,j)$. 
  Therefore the assertion follows. 

  \par
  (2) 
  By Lemma \ref{claim:restriction}, 
  we may prove that if $\# V = i+j$ and $\beta_{i-1, i+j} (I(G)) \neq 0$, 
  then $\beta_{i-1, i+j} (I(G)) = 1$. 
  We proceed the proof by induction on $\# E(G)$. 
  When $\# E(G) = 1$, the claim is trivial. We assume that $\# E(G) \geq 2$. 
  Note that the equality (\ref{eq:HaTuylproof2}) holds also in this case. 
  Since $G$ is a forest, the neighbourhood of $v$ in $G$ consists of only $u$. 
  Thus $v$ is an isolated vertices of $G \setminus e$.  
  This implies that 
  $\beta_{i-1, i+j} (I(G \setminus e)) = 0$. 
  Thus $\beta_{i-2-t, (i-1-t)+(j-1)} (I(G')) \neq 0$ 
  by (\ref{eq:HaTuylproof2}). 
  Then we have $\beta_{i-2-t, (i-1-t)+(j-1)} (I(G')) = 1$ 
  by inductive hypothesis and we obtain  
  $\beta_{i-1, i+j} (I(G)) = 1$ again by (\ref{eq:HaTuylproof2}), 
  as desired. 
\end{proof}

\par
The next example shows that Theorem \ref{claim:Betti_chordal} (1) is false 
for a general graph. 
\begin{exam}
  Let us consider the following non-chordal graph $G_1$ on the vertex set 
  $V = \{ 1, 2, \ldots, 6 \}$: 
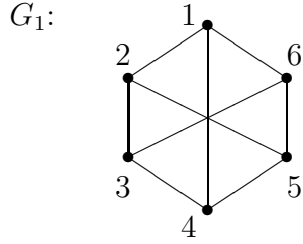
  \begin{center}
  \begin{picture}(115,110)(-25,0)
    \put(-45,90){$G_1$:}
    \put(0,70){\circle*{4}}
    \put(0,40){\circle*{4}}
    \put(30,20){\circle*{4}}
    \put(30,90){\circle*{4}}
    \put(60,70){\circle*{4}}
    \put(60,40){\circle*{4}}
    \put(-5,75){$2$}
    \put(60,75){$6$}
    \put(-5,25){$3$}
    \put(60,25){$5$}
    \put(20,10){$4$}
    \put(20,90){$1$}
    \put(30,20){\line(-3,2){30}}
    \put(30,20){\line(3,2){30}}
    \put(30,90){\line(-3,-2){30}}
    \put(30,90){\line(3,-2){30}}
    \put(0,40){\line(0,1){30}}
    \put(60,40){\line(0,1){30}}
    \put(30,20){\line(0,1){70}}
    \put(0,40){\line(2,1){60}}
    \put(0,70){\line(2,-1){60}}
  \end{picture}
  \end{center}
  Actually, $G_1$ is the complete bipartite graph $K_{3,3}$. In particular 
  it is unmixed. The Betti diagram of $S/I(G_1)$ is 
  \begin{center}
  \begin{tabular}{c|rrrrrr}
    $j \backslash i$ & $0$ & $1$ & $2$ & $3$ & $4$ & $5$ \\ \hline
    $0$ & $1$ & & & & & \\
    $1$ & & $9$ & $18$ & $15$ & $6$ & $1$ 
  \end{tabular}
  \end{center}
  Here the entry of $j$th row and $i$th column stands for 
  $\beta_{i,i+j} (S/I(G_1))$. 
  There is no bouquet with $4$ flowers in $G_1$ though 
  $\beta_{4,4+1} (S/I(G_1)) \neq 0$. Also $\pd S/I(G_1) = 5$ 
  while $d_{G_1} = 3$. 
\end{exam}

\par
The next example shows that Theorem \ref{claim:Betti_chordal} (2) is false 
for chordal graphs. 
\begin{exam}
  Let $G_2$ be the chordal graph with the vertex set $V = \{ 1, 2, 3 \}$ 
  and the edge set $\{ \{ 1,2 \}, \{ 1,3 \}, \{ 2,3 \} \}$: 
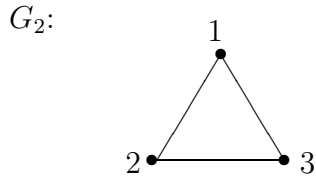
  \begin{center}
  \begin{picture}(100,90)
    \put(-20,70){$G_2$:}
    \put(34,20){\circle*{4}}
    \put(84,20){\circle*{4}}
    \put(60,60){\circle*{4}}
    \put(24,15){$2$}
    \put(90,15){$3$}
    \put(55,65){$1$}
    \put(34,20){\line(1,0){48}}
    \put(60,60){\line(3,-5){24}}
    \put(60,60){\line(-3,-5){24}}
  \end{picture}
  \end{center}
  The Betti diagram of $S/I(G_2)$ is 
  \begin{center}
  \begin{tabular}{c|rrr}
    $j \backslash i$ & $0$ & $1$ & $2$ \\ \hline
    $0$ & $1$ & & \\
    $1$ & & $3$ & $2$ 
  \end{tabular}
  \end{center}
  Although $\beta_{2,2+1} (S/I(G_2)) = 2$, a subset of $V$ with cardinality 
  $2+1$ is only $V$. 
\end{exam}

\section{Another characterization of the projective dimension}
We gave a characterization of the projective dimension for chordal graphs 
in Theorem \ref{claim:Betti_chordal}. However it seems that it is not easy 
to read it from the graph. In this section, we introduce another notion 
of disjointness of a set of bouquets, 
which is rather weaker than strongly disjointness and give 
another characterization of the projective dimension for chordal graphs 
using this notion. 

\par
Let $G$ be a finite simple graph and 
$\mathcal{B} = \{ B_1, B_2, \ldots, B_j \}$ a set of bouquets of $G$. 
\begin{defn}
  \label{defn:SemiStrongDisjoint}
  We say that a set $\mathcal{B} = \{ B_1, B_2, \ldots, B_j \}$ of 
  bouquets of $G$ is \textit{semi-strongly disjoint} in $G$ 
  if the following $2$ conditions are satisfied: 
  \begin{enumerate}
  \item $V(B_k) \cap V(B_{\ell}) = \emptyset$ for all $k \neq \ell$. 
  \item Any two vertices belonging to $R({\mathcal{B}})$ 
    are not adjacent in $G$. 
  \end{enumerate}
\end{defn}

\par
As noted in Remark \ref{rmk:root}, 
if $\mathcal{B}$ is strongly disjoint, then 
$\mathcal{B}$ is also semi-strongly disjoint. 
%
We set 
\begin{displaymath}
  d_G' := \max \{ \# F({\mathcal{B}}) \; : \; 
  \text{$\mathcal{B}$ is a semi-strongly disjoint set of bouquets of $G$} \}. 
\end{displaymath}
In general, the inequality $d_G \leq d_G'$ holds. 
There exists a graph $G$ with $d_G < d_G'$ as the following example shows. 
\begin{exam}
  \label{eq:d<d'}
  Let $G_3$ be the following graph: 
  
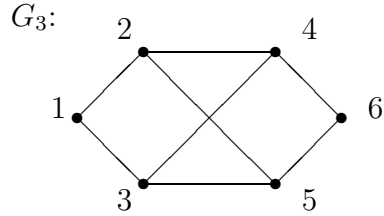
\begin{center}
  \begin{picture}(150,80)
    \put(0,70){$G_3$:}
    \put(50,60){\circle*{4}}
    \put(100,60){\circle*{4}}
    \put(50,10){\circle*{4}}
    \put(100,10){\circle*{4}}
    \put(25,35){\circle*{4}}
    \put(125,35){\circle*{4}}
    \put(40,65){$2$}
    \put(110,65){$4$}
    \put(40,0){$3$}
    \put(110,0){$5$}
    \put(15,35){$1$}
    \put(135,35){$6$}
    \put(25,35){\line(1,1){25}}
    \put(25,35){\line(1,-1){25}}
    \put(125,35){\line(-1,1){25}}
    \put(125,35){\line(-1,-1){25}}
    \put(50,60){\line(1,0){50}}
    \put(50,60){\line(1,-1){50}}
    \put(50,10){\line(1,0){50}}
    \put(50,10){\line(1,1){50}}
  \end{picture}
  \end{center}
  Note that $G_3$ is a bipartite graph which is not unmixed. 
  It is easy to see that the distance of any two edge of $G_3$ is at most $2$. 
  That is there exists no $3$-disjoint edges in $G_3$. 
  Thus the strongly disjoint set of bouquets of $G_3$ 
  consists of only one bouquet. In particular $d_{G_3} = 3$. 
  On the other hand let $B_1$ (resp.\  $B_2$) be the bouquet whose 
  root is $1$ (resp.\  $6$) and whose flowers are $2,3$ (resp.\  $4,5$). 
  Since $\{ 1, 6 \}$ is not an edge of $G_3$, the set of bouquets $B_1, B_2$ 
  is semi-strongly disjoint and $d_{G_3}' = 4 > 3 = d_{G_3}$. 
\end{exam}

The following result is the main theorem of this section. 
\begin{thm}
  \label{claim:pd}
  Let $G$ be a chordal graph.  Then 
  \begin{displaymath}
    \pd S/I(G) = d_G = d_G'. 
  \end{displaymath}
\end{thm}
\begin{proof}
  By Theorem \ref{claim:Betti_chordal}, we may prove that $d_G \geq d_G'$. 
  We use induction on $\# E(G)$. When $\# E(G) = 1$, clearly 
  $d_G = d_G' = 1$ hold. 

  \par
  We consider the case of $\# E(G) \geq 2$. 
  Let $\mathcal{B} = \{ B_1, B_2, \ldots, B_j \}$ be a semi-strongly 
  disjoint set of bouquets of $G$ with $\# F({\mathcal{B}}) = d_G'$. 
  When $j=1$, the assertion is clear because in this case, 
  $\mathcal{B}$ is also strongly disjoint. We assume that $j \geq 2$. 
  Let $e = \{ u, v \}$ be an edge of $G$ such that
  $G_{N(v)}$ is a complete graph as in Lemma \ref{claim:HaTuyl}. 
  If $e \notin S({\mathcal{B}})$, then $\mathcal{B}$ is a set of bouquets of 
  $G \setminus e$ and semi-strongly disjoint also in $G \setminus e$. 
  Hence by inductive hypothesis and Lemma \ref{claim:d} (1), we have 
  \begin{displaymath}
    d_G' \leq d_{G \setminus e}' = d_{G \setminus e} \leq d_G, 
  \end{displaymath}
  as desired. 

  \par
  Next we consider the case where $e \in S({\mathcal{B}})$, say $e \in E(B_j)$. 
  First we consider the case where $u$ is a neighbour of the root of $B_1$. 
  In this case the root of $B_j$ must be $v$. 
  Let $B_1'$ be the bouquet obtained by adding $u$ to 
  $B_1$ as a flower and $B_j'$ the one obtained by removing $u$ from $B_j$. 
  Then $\mathcal{B}' := \{ B_1', B_2, \ldots, B_{j-1}, B_{j}' \}$ 
  is a semi-strongly disjoint set of bouquets of $G$ 
  with $e \notin S({\mathcal{B}})$. 
  Thus we have $d_G \geq d_G'$ as shown in the above. 
  Therefore we may assume that $u$ is not a neighbour of the root $w_k$ 
  of $B_k$ for any $k=1, 2, \ldots, j-1$. 
  Then since $G_{N(v)}$ is a complete graph, 
  we may assume that $u$ is a root of $B_j$. 
  Let $\mathcal{B}'$ be the set of bouquets which is obtained by 
  removing all flowers those are neighbours of $u$ 
  from the set of bouquets $B_1, B_2, \ldots, B_{j-1}$. 
  Then $\mathcal{B}'$ is a semi-strongly disjoint set of bouquets of 
  $G' = G_{V \setminus (N(u) \cup \{ u \})}$ and 
  $d_G' - (t+1) \leq \# F({\mathcal{B}}') \leq d_{G'}'$. 
  Combining this with the inductive hypothesis and Lemma \ref{claim:d} (2), 
  we have 
  \begin{displaymath}
    d_G' \leq d_{G'}' + (t+1) = d_{G'} + (t+1) \leq d_G. 
  \end{displaymath}
  This completes the proof. 
\end{proof}

\par
We cannot replace the strongly disjointness by the semi-strongly one 
on a characterization of the graded Betti numbers 
for a forest in Theorem \ref{claim:Betti_chordal} (2). 
\begin{exam}
  Let $G_4$ be the line graph with $7$ vertices and $\mathcal{B}$ 
  the set of bouquets of $G_4$ as the following picture: 
  
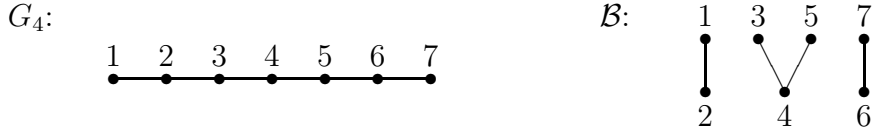
\begin{center}
    \begin{picture}(230,70)(0,35)
      \put(10,80){$G_{4}$:}
      \put(50,60){\line(1,0){120}}
      \multiput(50,60)(20,0){7}{\circle*{4}}
      \put(47,65){$1$}
      \put(67,65){$2$}
      \put(87,65){$3$}
      \put(107,65){$4$}
      \put(127,65){$5$}
      \put(147,65){$6$}
      \put(167,65){$7$}
    \end{picture}
    \begin{picture}(120,70)(-50,20)
      \put(-50,65){$\mathcal{B}$:}
      \put(-10,40){\circle*{4}}
      \put(-13,27){$2$}
      \put(-10,60){\circle*{4}}
      \put(-13,65){$1$}
      \put(-10,40){\line(0,1){20}}
      \put(20,40){\circle*{4}}
      \put(17,27){$4$}
      \put(10,60){\circle*{4}}
      \put(7,65){$3$}
      \put(30,60){\circle*{4}}
      \put(27,65){$5$}
      \put(20,40){\line(-1,2){10}}
      \put(20,40){\line(1,2){10}}
      \put(50,40){\circle*{4}}
      \put(47,27){$6$}
      \put(50,60){\circle*{4}}
      \put(47,65){$7$}
      \put(50,40){\line(0,1){20}}
    \end{picture}
  \end{center}
  It is easy to see that $\mathcal{B}$ is semi-strongly disjoint. 
  But $\mathcal{B}$ is not strongly-disjoint. 
  To see this, let $B$ be the bouquet of $\mathcal{B}$ whose root is $4$. 
  Then the stem of $B$ which is $3$-disjoint with $\{ 1,2 \}$ in $G_4$ is only 
  $\{ 4,5 \}$ but it is not $3$-disjoint with $\{ 6,7 \}$ in $G_4$. 

  \par
  In fact, the type of $\mathcal{B}$ is $(4,3)$ while the Betti diagram of 
  $S/I(G_4)$ is given by 
  \begin{center}
  \begin{tabular}{c|rrrrr}
    $j \backslash i$ & $0$ & $1$ & $2$ & $3$ & $4$ \\ \hline
    $0$ & $1$ & & & & \\
    $1$ & & $6$ & $5$ & & \\
    $2$ & & & $6$ & $9$ & $3$ 
  \end{tabular}
  \end{center}
  and $\beta_{4, 4+3} (S/I(G_4)) = 0$. 
\end{exam}

\begin{remark}
In general, inequalities $\pd S/I(G) \geq d_G$ 
(Theorem \ref{claim:NonVanishBetti}) and $d_G' \geq d_G$ holds. 
Then are there any relations between $\pd S/I(G)$ and $d_G'$? 
The inequality $\pd S/I(G) \geq d_G'$ holds for many graphs $G$. 
But we do not know whether this is always true or not. 
\end{remark}

\par
Finally, we give an application of Theorem \ref{claim:pd}. 
Let $G$ be a finite simple graph on $V$. 
A subset $\mathcal{C}$ is said to be a \textit{vertex cover} of $G$
if it intersects all edges of $G$. 
We say that a vertex cover of $G$ is \textit{minimal} 
if it has no proper subset which is also a vertex cover of $G$. 
There is one-to-one correspondence between the vertex covers of $G$ 
and the minimal prime divisors of $I(G)$: 
\begin{displaymath}
  I(G) = \bigcap_{\text{$\mathcal{C}$: a minimal vertex cover of $G$}} 
    (x_i \; : \; i \in \mathcal{C}). 
\end{displaymath}
The \textit{big height} of $I(G)$, denoted by $\bight I(G)$, 
is defined as the maximum heights 
among the minimal prime divisors of $I(G)$. 
In general, the inequality $\pd S/I(G) \geq \bight I(G)$ holds. 
As a corollary of Theorem \ref{claim:pd}, we recover the following result 
which follows from Morey and Villarreal \cite[Corollary 3.33]{MV} 
together with Francisco and Van Tuyl \cite[Theorem 3.2]{FranTuyl}; 
see also Herzog, Hibi, and Zheng \cite{HHZ}. 
\begin{cor}[{\cite{MV}, \cite{FranTuyl}, \cite{HHZ}}]
  \label{claim:bight}
  Let $G$ be a chordal graph on $V$. We set $S = K[V]$. 
  Then 
  \begin{displaymath}
    \bight I(G) = \pd S/I(G). 
  \end{displaymath}
  In particular, if $I(G)$ is unmixed, then $S/I(G)$ is Cohen--Macaulay. 
\end{cor}


\begin{proof}
  We may prove that $\pd S/I(G) \leq \bight I(G)$. 
  Let $\mathcal{B} = \{ B_1, B_2, \ldots, B_j \}$ be 
  a semi-strongly disjoint set of bouquets of $G$ with 
  $\# F ({\mathcal{B}}) = d_G'$. 
  Then $\# F({\mathcal{B}}) = \pd S/I(G)$ by Theorem \ref{claim:pd}. 
  We claim that $F({\mathcal{B}})$ is a minimal vertex cover of $G$. 
  If $F({\mathcal{B}})$ is a vertex cover of $G$, then the minimality 
  of it is clear. Hence we only show that $F({\mathcal{B}})$ is 
  a vertex cover of $G$. 

  \par
  Suppose that, on the contrary, there exists an edge $e = \{ u,v \} \in E(G)$ 
  with $e \cap F({\mathcal{B}}) = \emptyset$. 
  When $u$ is a neighbour of the root of $B_k \in \mathcal{B}$ for some $k$, 
  the set $\mathcal{B}'$ of bouquets obtained by adding $u$ to $B_k$ 
  as a flower is semi-strongly disjoint 
  and satisfies $\# F(\mathcal{B}') = d_G'+1$. 
  This is a contradiction. 
  Therefore $u, v \notin F({\mathcal{B}}) \cup R({\mathcal{B}})$. 
  Let $B_e$ denotes the bouquet of $G$ with only one stem $e$. 
  Then $\mathcal{B}'' := \mathcal{B} \cup \{ B_e \}$ is a 
  semi-strongly disjoint set of bouquets of $G$ and 
  $\# F(\mathcal{B}'') = d_G' +1$. This is also a contradiction. 
\end{proof}

\begin{acknowledgement}
  The author thanks Professors J\"{u}rgen Herzog and Takayuki Hibi 
  for the suggestion to consider the characterization 
  of the graded Betti numbers. 
  She also thanks Professor Akimichi Takemura 
  for giving her attention to chordal graphs. 
  The Author is grateful to Professor Naoki Terai for telling 
  her Corollary \ref{claim:bight}. 
  She is also grateful to Professor Rafael H. Villarreal for pointing out 
  that Corollary \ref{claim:bight} was proved in a more general setting in 
  \cite{MV}. The author thanks Professor Russ Woodroofe for giving her useful 
  suggestions and comments. 
\end{acknowledgement}


\end{document}